\documentstyle[amssymb,12pt]{amsart}
\newtheorem{prop}{Proposition}[section]
\newtheorem{prop:def}{Proposition-Definition}[section]
\newtheorem{defin}{Definition}[section]

\newtheorem{thm}{Theorem}[section]

\theoremstyle{remark}
\newtheorem{remark}{Remark}
\newtheorem{example}{Example}
\newtheorem{notation}{Notation}
\newtheorem{note}{Note}

\begin{document}

\newcommand{\nc}{\newcommand} \nc{\on}{\operatorname}

\nc{\pa}{\partial}

\nc{\cA}{{\cal A}} \nc{\cB}{{\cal B}}\nc{\cC}{{\cal C}} 
\nc{\cE}{{\cal E}} \nc{\cG}{{\cal G}}\nc{\cH}{{\cal H}} 
\nc{\cI}{{\cal I}} \nc{\cJ}{{\cal J}}\nc{\cK}{{\cal K}} 
\nc{\cL}{{\cal L}} \nc{\cR}{{\cal R}} \nc{\cS}{{\cal S}}   
\nc{\cV}{{\cal V}}  \nc{\cX}{{\cal X}}

\nc{\sh}{\on{sh}}\nc{\Id}{\on{Id}}\nc{\id}{\on{id}}
\nc{\Diff}{\on{Diff}}
\nc{\Perm}{\on{Perm}}\nc{\conc}{\on{conc}}\nc{\Alt}{\on{Alt}}
\nc{\ad}{\on{ad}}\nc{\Der}{\on{Der}}\nc{\End}{\on{End}}
\nc{\no}{\on{no\ }} \nc{\res}{\on{res}}\nc{\ddiv}{\on{div}}
\nc{\Sh}{\on{Sh}} \nc{\card}{\on{card}}\nc{\dimm}{\on{dim}}
\nc{\Sym}{\on{Sym}} \nc{\Jac}{\on{Jac}}\nc{\Ker}{\on{Ker}}
\nc{\Vect}{\on{Vect}} \nc{\Spec}{\on{Spec}}\nc{\Cl}{\on{Cl}}
\nc{\Imm}{\on{Im}}\nc{\limm}{\lim}\nc{\Ad}{\on{Ad}}
\nc{\ev}{\on{ev}} \nc{\Hol}{\on{Hol}}\nc{\Det}{\on{Det}}
\nc{\Bun}{\on{Bun}}\nc{\diag}{\on{diag}}\nc{\pr}{\on{pr}} 
\nc{\Span}{\on{Span}}\nc{\Comp}{\on{Comp}}\nc{\Part}{\on{Part}}
\nc{\tensor}{\on{tensor}}\nc{\ind}{\on{ind}}\nc{\Tors}{\on{Tors}}

\nc{\al}{\alpha}\nc{\g}{\gamma}\nc{\de}{\delta}
\nc{\eps}{\epsilon}\nc{\la}{{\lambda}}
\nc{\si}{\sigma}\nc{\z}{\zeta}

\nc{\La}{\Lambda}

\nc{\ve}{\varepsilon} \nc{\vp}{\varphi} 

\nc{\AAA}{{\mathbb A}}\nc{\CC}{{\mathbb C}}\nc{\ZZ}{{\mathbb Z}} 
\nc{\QQ}{{\mathbb Q}} \nc{\NN}{{\mathbb N}}\nc{\VV}{{\mathbb V}} 
\nc{\KK}{{\mathbb K}} 

\nc{\ff}{{\mathbf f}}\nc{\bg}{{\mathbf g}}
\nc{\ii}{{\mathbf i}}\nc{\kk}{{\mathbf k}}
\nc{\bl}{{\mathbf l}}\nc{\zz}{{\mathbf z}} 
\nc{\pp}{{\mathbf p}}\nc{\qq}{{\mathbf q}}

\nc{\cF}{{\cal F}}\nc{\cM}{{\cal M}}\nc{\cO}{{\cal O}}
\nc{\cT}{{\cal T}}\nc{\cW}{{\cal W}}

\nc{\ub}{{\underline{b}}}
\nc{\uk}{{\underline{k}}} \nc{\ul}{{\underline}}
\nc{\un}{{\underline{n}}} \nc{\um}{{\underline{m}}}
\nc{\up}{{\underline{p}}}\nc{\uq}{{\underline{q}}}
\nc{\us}{{\underline{s}}}\nc{\ut}{{\underline{t}}}
\nc{\uw}{{\underline{w}}}
\nc{\uz}{{\underline{z}}}
\nc{\ual}{{\underline{\alpha}}}\nc{\ualpha}{{\underline{\alpha}}}
\nc{\ugamma}{{\underline{\gamma}}}
\nc{\ula}{{\underline{\lambda}}}\nc{\umu}{{\underline{\mu}}}
\nc{\unu}{{\underline{\nu}}}\nc{\usigma}{{\underline{\sigma}}}
\nc{\utau}{{\underline{\tau}}}
\nc{\uN}{{\underline{N}}}\nc{\uM}{{\underline{M}}}
\nc{\uK}{{\underline{K}}}

\nc{\A}{{\mathfrak a}} \nc{\B}{{\mathfrak b}} \nc{\G}{{\mathfrak g}}
\nc{\D}{{\mathfrak d}} \nc{\HH}{{\mathfrak h}}  \nc{\iii}{{\mathfrak
i}}   \nc{\mm}{{\mathfrak m}} \nc{\N}{{\mathfrak
n}}\nc{\ttt}{{\mathfrak{t}}}  \nc{\U}{{\mathfrak u}}\nc{\V}{{\mathfrak
v}}

\nc{\SL}{{\mathfrak{sl}}}

\nc{\SG}{{\mathfrak S}}

\nc{\wt}{\widetilde} \nc{\wh}{\widehat}
\nc{\bn}{\begin{equation}}\nc{\en}{\end{equation}} \nc{\td}{\tilde}

%
%
%

\newcommand{\ldar}[1]{\begin{picture}(10,50)(-5,-25)
\put(0,25){\vector(0,-1){50}}
\put(5,0){\mbox{$#1$}} 
\end{picture}}

\newcommand{\luar}[1]{\begin{picture}(10,50)(-5,-25)
\put(0,-25){\vector(0,1){50}}
\put(5,0){\mbox{$#1$}}
\end{picture}}

\title[]{Quantum torsors}
\author{Cyril Grunspan}
\address{UFR de Math\'ematique et Informatique, Universit\'e de Strasbourg, IRMA,
7 rue Ren\'e Descartes, F-67084 Strasbourg cedex (France)}
\email{grunspan@@math.u-strasbg.fr}

\begin{abstract}
  The following text is a short version of a forthcoming preprint about
  torsors. The adopted viewpoint is an old reformulation of torsors recalled
  recently by Kontsevich {\cite{Kon}}. We propose a unification of the
  definitions of torsors in algebraic geometry and in Poisson manifolds
  (Example \ref{geoalg} and section \ref{Poissontor}). We introduce the notion
  of a quantum torsor (Definition \ref{tor}). Any quantum torsor is equipped
  with two comodule-algebra structures over Hopf algebras and these structures
  commute with each other (Theorem \ref{premtheo}.) In the finite dimensional
  case, these two Hopf algebras share the same finite dimension (Proposition
  \ref{hopfgalprop}). We show that any Galois extension of a field is a torsor
  (Example \ref{abel}) and that any torsor is a Hopf-Galois extension (section
  \ref{toto}). We give also examples of non-commutative torsors without
  character (Example \ref{examp}). Torsors can be composed (Theorem
  \ref{zorro}). This leads us to define for any Hopf algebra, a new
  group-invariant, its torsors invariant (Theorem \ref{invtor}). We show how
  Parmentier's quantization formalism of ``affine Poisson groups'' is part of
  our theory of torsors (Theorem \ref{parmi}). 
\end{abstract}
\thanks{I am very grateful to B. Enriquez for
his help and his support. I am also grateful to P. Cartier, A. Chambert-Loir,
S. Natale and S. Parmentier for comments and remarks.}
\maketitle


\section{Introduction}

\subsection{General overview}

The aim of our work is to give a meaning and to develop a general theory for
quantum torsors, starting from the principle that most of the objects of
``traditional'' commutative geometry should have a counterpart in the
framework of non-commutative geometry. In algebraic geometry, torsors are
familiar objects. Indeed, we know that they are linked with the problem of
inner forms of algebraic groups and classified with the help of Galois
cohomology groups : modulo an equivalence relation, torsors over an algebraic
group $G$ defined over a field $k$ are in correspondence with a pointed
cohomology set $H^1 ( \text{Gal} ( \bar{k} / k ), G ( \bar{k} ))$ where $k$ is an
algebraic closure of $k$, such that the trivial torsors $X$ on $k$ (the ones
with $X ( k ) \not=\emptyset )$ correspond to the trivial cocyle.
Therefore, it seems natural to develop a similar theory in the framework of
quantum groups or Hopf algebras. On the other hand, in the category of Poisson
manifolds, torsors are well-known objects and have been classified by Dazord
and Sondaz {\cite{DS}} under the name of ``affine Poisson groups'' (for a
precursor, see {\cite{STS}}) and quantized by Parmentier {\cite{Par}}.
However, the two definitions of torsors in algebraic geometry and in Poisson
geometry do not coincide. In algebraic geometry, a torsor is a scheme $X$
equipped with a group-action $m : G \times X \longrightarrow X$ where $G$ is a
group-scheme such that the map $( m \times \text{pr}_X ) : G \times X
\longrightarrow X \times X$ is an isomorphism of schemes, whereas in the
category of Poisson manifolds, a torsor $X$ is (still) the data of a
Poisson-Lie group $G$ and a group-action $m : G \times X \longrightarrow X$
plus certain conditions saying that the stabilizers of the action are all
trivial, but in that case, the bijective map $( m \times \text{pr}_X ) : G
\times X \longrightarrow X \times X$ is {\em{not}} a Poisson map.

Thus, we are looking for unifying in a same theory the results of algebraic
geometry about torsors and Parmentier's quantization results. Our starting
point is an old intrinsic reformulation of affine structures originally
suggested by Baer {\cite{B}} and developed later on by Certaine {\cite{C}},
Vagner {\cite{Vag}}, Kock {\cite{Koc}}, Weinstein {\cite{W}} and recalled
recently by Kontsevich {\cite{Kon}}. One of the corollaries of our work is the
definition of a non trivial group-invariant $\text{Tor} ( H )$ associated with
any Hopf algebra $H$ (Theorem \ref{invtor}).

\subsection{\label{secd}Definition of a torsor}

In the category of sets, affine algebraic manifolds over an algebraically
closed field or Poisson manifolds, a torsor is a $G$-{\em{principal
homogeneous space}} where $G$ is a group in the category. Even though the
axioms asserting that $G$ acts transitively can be quantized without any
problem, the axiom saying that all stabilizers are trivial is not so easy to
quantize, especially when the quantum manifolds have no point (i.e., when the
algebras of functions have no character). According to the point of view
originally developed by Baer {\cite{B}}, a (classical) torsor is the data of
an object $X$ and a composition law $\mu_X : X^3 \longrightarrow X$ satisfying
some associativity relations, called {\em{parallelogram}} relations and
analogous to the ones we would get by taking $X = G$ a group and
\begin{equation}
  \begin{array}{cccc}
    \mu_G : & G^3 & \longrightarrow & G\\
    & ( g, h, k ) & \longmapsto & g h^{- 1} k
  \end{array} \label{defmug}
\end{equation}
By reversing the arrows, Kontsevich noted that it was possible to define in
this way all the torsors of the algebraic geometry {\cite{Kon}}. We will see
that it is also possible to extend this idea in non-commutative geometry.

\subsection{Classical torsors}

In pointed categories such as the category of sets, affine algebraic manifolds
over an algebraically closed field or Poisson manifolds, it is easy to see
that the map $\mu_X : X^3 \longrightarrow X$ should verify the following
relations :
\begin{eqnarray}
  \forall a, b, c, d, e \in X,\quad \mu_X ( a, a, b ) & = & b \\
  \mu_X ( a, b, b ) & = & a \\
  \mu_X ( \mu_X ( a, b, c ), d, e ) & = & \mu_X ( a, b, \mu_X ( c, d, e ))\\
  & = & \mu_X ( a, \mu_X ( d, c, b ), e )
\end{eqnarray}
Conversely, if an object $X$ is equipped with a law $\mu_X$ satisfying these
equalities, then we can easily find two groups of the category $G_l ( X )$ and
$G_r ( X )$ acting simply and transitively on $X$, the first one from the left
and the second one from the right and $X$ is a {\em{principal homogeneous
space}} on $G_l ( X )$ or $G_r ( X )$. Furthermore, we are able to classify
all the classical torsors. In the category of sets, as well as in the category
of affine algebraic manifolds over an algebraically closed field, any torsor
is isomorphic to the trivial torsor $( G, \mu_G )$ on a group $G$. In the
category of Poisson manifolds, it can be shown that any torsor is isomorphic
to an ``affine Poisson group'' according to the terminology introduced by
Dazord and Sondaz {\cite{DS}} and thus is identified with a triple $(
\mathfrak{g}, \delta, f )$ where $( \mathfrak{g}, \delta )$ is a Lie bialgebra
and $f \in \Lambda^2 ( \mathfrak{g} )$ is a ``classical Drinfeld twist'' for
$\delta$. The affine Poisson group associated with this triple is defined in
the following way : denoting by $( G, P_G )$ the connected, simply connected
Poisson-Lie group associated with the bialgebra $( \mathfrak{g}, \delta )$,
the affine Poisson group is $G$ as a manifold equipped with the
{\em{affine}} Poisson structure given by the Poisson bivector $\pi :=
P_G + f^L$ where $f^L$ is the left translation of $f$ on $G$. On this Poisson
manifold, the Poisson-Lie groups acting simply and transitively are $( G, P_G
)$ by translation from the left and $( G, P_G + f^L - f^R )$ by translation
from the right ($f^R$ being the right translation of $f$ on $G$).

\section{Quantum torsors}

\subsection{Non-commutative torsors}

Let us consider now a commutative field $k$ and a commutative $k$-algebra $A$
without any zero divisor. For example, $A = k$ or $A = k [[ \hbar ]]$. We will
work in the category of $A$-(unitary associative) algebras.

\begin{defin}
  \label{tor}An $A$-torsor is a quintuple $( T, m_T, 1_T, \mu_T, \theta_T )$
  where $( T, m_T, 1_T )$ is an $A$-algebra, $\mu_T : T \longrightarrow T
  \otimes_A T^{\text{op}} \otimes_A T$ is an $A$-algebra morphism and
  $\theta_T : T \longrightarrow T$ is an $A$-algebra automorphism satisfying
  the following axioms :
  \begin{eqnarray}
    \forall x \in T, \quad ( \text{Id}_T \otimes m_T ) \circ \mu_T ( x ) & = &
    x \otimes 1_T \\
    ( m_T \otimes \text{Id}_T ) \circ \mu_T ( x ) & = & 1_T \otimes x\\
    ( \text{Id}_T \otimes \text{Id}_{T^{\text{op}}} \otimes \mu_T ) \circ
    \mu_T & = & ( \mu_T \otimes \text{Id}_{T^{\text{op}}} \otimes \text{Id}_T
    ) \circ \mu_T\\
    \theta^{(3)}_T\circ ( \mu_T \otimes
    \text{Id}_{T^{\text{op}}} \otimes \text{Id}_T ) \circ \mu_T & = & (
    \text{Id}_T \otimes \mu_T^{\text{op}} \otimes \text{Id}_T ) \circ \mu_T
    \label{theta}\\
    ( \theta_T \otimes \theta_T \otimes \theta_T ) \circ \mu_T & = & \mu_T
    \circ \theta_T 
  \end{eqnarray}
  with $\mu_T^{\text{op}} := \tau_{( 13 )} \circ \mu_T$ and 
  $\theta^{(3)}_T:= ( \text{Id}_T \otimes \text{Id}_{T^{\text{op}}} \otimes \theta_T \otimes
    \text{Id}_{T^{\text{op}}} \otimes \text{Id}_T )$. If $m_T =
  m_T^{\text{op}}$, the torsor is said to be commutative. If $\mu_T =
  \mu_T^{\text{op}}$, the torsor is said to be endowed with a commutative law.
\end{defin}

\begin{note}
  If $( T, m_T, 1_T, \mu_T, \theta_T )$ is an $A$-torsor, then $\theta_T$ is
  fully determined by $m_T$ and $\mu_T$. For instance, if the torsor is either
  commutative or endowed with a commutative law, then $\theta_T = \text{Id}_T.$
\end{note}  

  \begin{remark}
    Of course, a given algebra needs not carry a torsor structure. For
    example, if $\text{char} ( k ) \not= 2$, then there is no $k$-torsor
    structure on the $k$-algebra $k [ X ] / ( X^2 )$.
  \end{remark}  
 
  \begin{remark}
      \label{opp}If $(T, m_T, 1_T, \mu_T, \theta_T )$ is an $A$-torsor, then
      $( T^{\text{op}}, m_T^{\text{op}}, 1_T, \mu_T^{\text{op}}, \theta_T )$
      is also an $A$-torsor, called its opposite torsor. 
  \end{remark}

\begin{notation} 
We will use generalized Sweedler notations. If $( T, m_T, 1_T, \mu_T, \theta_T
)$ is an $A$-torsor, then for all $x \in T$, forgetting the symbol $\sum,$we
denote $\mu_T ( x ) = x^{( 1 )} \otimes x^{( 2 )} \otimes x^{( 3 )}$
and$\mu_T^{( n )} ( x ) = x^{( 1 )} \otimes \ldots \otimes x^{( 2 n + 1 )}, n
\in \mathbb{N}$, where $\mu_T^{( n )}$ satisfies the induction $\mu_T^{( 0 )}
:= \text{Id}_T$ and $\mu_T^{( n )} = ( \mu_T^{( n - 1 )} \otimes
\text{Id}_{T^{\text{op}}} \otimes \text{Id}_T ) \circ \mu_T$. 
\end{notation}

Then, the torsor axioms show that for any {\em{odd}} integer $i$, we have :
\[ x^{( 1 )} \otimes \ldots \otimes x^{( i - 1 )} \otimes x^{( i )^{( 1 )}}
   \otimes x^{( i )^{( 2 )}} \otimes x^{( i )^{( 3 )}} \otimes x^{( i + 1 )}
   \otimes \ldots \otimes x^{( 2 n - 1 )} = \mu_T^{( n )} ( x ) \]
and, for any {\em{even}} integer $i$,
\begin{eqnarray*}
  x^{( 1 )} \otimes \ldots \otimes x^{( i - 1 )} \otimes x^{( i )^{( 1 )}}
  \otimes x^{( i )^{( 2 )}} \otimes x^{( i )^{( 3 )}} \otimes x^{( i + 1 )}
  \otimes \ldots \otimes x^{( 2 n - 1 )} &  & \\
  = x^{( 1 )} \otimes \ldots \otimes x^{( i - 1 )} \otimes x^{( i + 2 )}
  \otimes \theta_T ( x^{( i + 1 )} ) \otimes x^{( i )} \otimes x^{( i + 3 )}
  \otimes \ldots \otimes x^{( 2 n + 1 )} &  & 
\end{eqnarray*}
In particular, we see that for all $x \in T, \theta_T ( x ) = x^{( 1 )} x^{(
2 )^{( 3 )}} x^{( 2 )^{( 2 )}} x^{( 2 )^{( 1 )}} x^{( 3 )}$.

As apparent with the study of Example \ref{tortriv} below which shows that
$\theta_T$ is an analogue of the square of the antipode in a Hopf algebra, it
seems necessary to introduce $\theta_T$ in the definition of a non-commutative
torsor.

{\bf{\begin{example}
  \label{tortriv}The trivial torsor of a Hopf algebra
\end{example}}}

Let $( H, m_H, \Delta_H, \eta_H, \varepsilon_H, S_H )$ be an $A$-Hopf algebra.
Then, $( H, m_H, 1_H, \mu_H, \theta_H )$ is an $A$-torsor with $1_H =
\eta_H ( 1 ), \mu_H := ( \text{Id}_H \otimes S_H \otimes \text{Id}_H )
\circ \Delta_H^{( 2 )}$ and $\theta_H := S_H^2 .$

{\bf{\begin{example}
  \label{geoalg}Affine torsors in algebraic geometry
\end{example}}}

As said in introduction, in algebraic geometry, a torsor is the data of a
scheme $X$, a group-scheme $G$ and an action $m : G \times X \longrightarrow
X$ such that the map :
\begin{equation}
  \text{ $( m \times \text{pr}_X ) : \quad G \times X \longrightarrow X \times
  X$} \label{deftorgeoalg}
\end{equation}
is an {\em{isomorphism of schemes}}. As noted by Kontsevich {\cite{Kon}},
in the language of {\em{affine}} scheme, $X$ correspond to a
comodule-algebra $A$ over an Hopf algebra $H$ and the map $\mu_X : X^3
\longrightarrow X$ obtained by composition :
\begin{equation}
  X^3 = X^2 \times X \tilde{\longrightarrow} G \times X \times X
  \longrightarrow G \times X \longrightarrow X
\end{equation}
where the second map is obtained by forgetting the second factor in $G \times
X \times X$ gives rise to a torsor-structure on $A$ (with of course $\theta_A
= \text{Id}_A )$. In other terms, the torsors of algebraic geometry are the
{\em{commutative}} torsors for the definition \ref{tor}. A good example of
torsor without any point that we should keep in mind is $X = \text{Isom} ( M (
n ), A )$ where $A$ is a {\em{simple central}} algebra. In this case, the
left group  (i.e., the group acting simply transitively from the left on the
torsor) is $\text{PGL} ( n )$ and the right group is $A^{\times} / k^{\times}$
if $k$ denotes the ground field.

{\bf{\begin{example}
  \label{quadra}Quadratic extensions of a field
\end{example}}}

The following examples of torsor structures on algebras (without character if
$d \not\in( k^{\times} )^2$) are nothing but {\em{subcases}} of the
previous example. If $k$ is a commutative field and if $d \in k^{\times},$then
the $k$-algebra $A = k [ X ] / ( X^2 - d )$ can be equipped with a $k$-torsor
structure by $\mu_A ( x ) = x \otimes x^{- 1} \otimes x, x = \text{cl} ( X )$
and $\theta_A = \text{Id}_A$. If $\text{char} ( k ) = 2$ and $A = k [ X ] / (
X^2 - d )$ or $A = k [ X ] / ( X^2 - X - d )$ with $d \in k$, the law $\mu_A$
defined on $A$ by $\mu_A ( x ) = 1 \otimes 1 \otimes x + 1 \otimes x \otimes 1
+ x \otimes 1 \otimes 1$ and $x = \text{cl} ( X )$ gives to $A$ a $k$-torsor
structure with $\theta_A = \text{Id}_A$.

{\bf{\begin{example}
  \label{abel}Generalization : Galois extension of a field
\end{example}}}

Let $K = k [ T ] / ( P )$ be a Galois extension of a field $k, t :=
\text{cl} ( T )$ a primitive element in $K,\, G := \text{Gal} ( K / k )$
the Galois group of the extension and $H := k^G := ( k G )^{\ast}$
the natural Hopf algebra of functions on $G$ with values in $k$. First, note
that there is a natural structure of $H^{\text{cop}}$-left comodule algebra on
$K$ given by a morphism $\Delta_K : K \longrightarrow H \otimes_k K$. Second,
note that we can identify $K \otimes_k K$ with $K [ T ] / ( P )$ such that
$t_1 := t \otimes 1 \in K \otimes_k K$ is identified to $\text{cl} ( T )
\in K [ T ] / ( P )$ and $t_2 := 1 \otimes t$ is identified to $t \in K
\subset K [ T ] / ( P )$. Third, note that there is a natural isomorphism
between $K^G := H \otimes_k K$ (the algebra of functions on $G$ with
values in $K$) and $K \otimes_k K \cong K [ T ] / ( P )$. Under this
isomorphism, $x \otimes y \in K \otimes_k K$ is associated to $\sum_{\sigma
\in G} 1_{\sigma} \otimes \sigma ( x ) y \in H \otimes K$ where $1_{\sigma}$
denotes the natural idempotent in $H$ associated to $\sigma$ by $< 1_{\sigma},
\tau > = \delta_{\sigma, \tau}$ for all $\tau \in G$ and $1_{\sigma} \otimes x
\in K^G$ is associated to $P_{\sigma} \times ( 1 \otimes x )$ with $P_{\sigma}
= \prod_{\tau \not{=} \sigma} \frac{t_1 - \tau ( t_2 )}{\sigma ( t_2 ) - \tau
( t_2 )}$ . Therefore, by composing $\Delta_K$ with the embedding of $H$ in
$K^G = H \otimes_k K \cong K \otimes_k K$, we get an algebra morphism :
\begin{equation}
  \begin{array}{rcl}
    \mu_K : K & \longrightarrow & K \otimes_k K \otimes_k K\\
    x & \longmapsto & \sum_{\sigma \in G} P_{\sigma} \otimes \sigma ( x )
  \end{array}
\end{equation}
We have $\mu_K^{\text{op}} ( x ) := \tau_{( 13 )} \circ \mu_K ( x ) =
\sum_{\sigma} \sigma ( x ) \otimes P_{\sigma^{- 1}}$ for all $x \in K$ and $(
g \otimes \text{Id} ) ( P_h ) = P_{h g^{- 1}}$ and $( \text{Id} \otimes g ) (
P_h ) = P_{g h}$ for any elements $g, h \in G$. From this, it can be shown
that $( K, m_K, 1_K, \mu_K, \text{Id}_K )$ is a $k$-torsor and that $G$ is a
group of automorphisms for the torsor $( K, m_K, 1_K, \mu_K, \text{Id}_K )$.
In others terms, Galois doesn't classify only torsors of the algebraic
geometry; Galois is also a torsor !

\begin{example}
  \label{examp} Non-commutative torsors without character
\end{example}

Let $k$ be a field and $n$ be a non-negative integer. Suppose that $k$
contains an element $q$ which is a $n$-th primitive root of $1$. For any
$\alpha$ and $\beta$ in $k^{\times}$, we denote by $A_{\alpha, \beta}$ the
{\em{non-commutative algebra}} with unit and {\em{without character}}
given by generators : $x, y$ and relations : $x^n = \alpha, y^n = \beta$ and
$x y = q y x$. The algebra $A_{\alpha, \beta}$ is a non-trivial {\em{cyclic
algebra}} and $\dim_k A_{\alpha, \beta} = n^2$. If $n = 2,$then $A_{\alpha,
\beta}$ is an algebra of quaternions. There is a natural structure of
non-commutative torsor endowed with a commutative law on $A_{\alpha, \beta}$
given by $\mu ( x ) = x \otimes x^{- 1}$$\otimes x,$$\mu ( y ) = y \otimes
y^{- 1} \otimes y$ and $\theta = \text{Id}$.

\subsection{\label{Poissontor}Poisson torsors}

\begin{defin}
  A $k$-Poisson torsor is a sextuple $( T, m_T, 1_T, \mu_T, \theta_T, \{, \}_T
  )$ such that 1) $( T, m_T, 1_T, \mu_T, \theta_T )$ is a commutative
  $k$-torsor, 2) $( T, m_T, 1_T, \{, \}_T )$ is a $k$-Poisson algebra and 3)
  the maps $\mu_T $and $\theta_T$ are Poisson maps where the Poisson structure
  on $T \otimes_k T\otimes_k T$ is given by the bracket :
  $\{ x \otimes y \otimes z, x' \otimes y' \otimes z' \} = \{ x, x'
     \} \otimes y y' \otimes z z' - x x' \otimes \{ y, y' \} \otimes z z' + x
     x' \otimes y y' \otimes \{ z, z' \}$
  for any elements $x,x',y,y',z,z'$ belonging to $T$.
\end{defin}

\begin{prop}
  If $G$ is a Poisson affine group, then $( G, \mu_G )$ is a classical Poisson
  torsor and $( \mathcal{O}_G, m_{\mathcal{O}_G}, 1_{\mathcal{O}_G} )$ can be
  naturally endowed with a Poisson torsor structure. Conversely, if $X$ is a
  smooth manifold and if $( \mathcal{O_{}}_X, m
  \mathcal{_{\mathcal{O_{}}_{}}}_{_X}, 1_{\mathcal{O}_X} )$ can be endowed
  with a Poisson torsor structure then $X$ is isomorphic to a Poisson affine
  group.
\end{prop}

However, Poisson torsor structures on commutative algebras without character
should exist (?), even if we are unable to give a single example. A way to get
such a torsor should be to consider an affine torsor $X$ without any point on
a field $k$ (cf Example \ref{geoalg}) and two $r$-matrices $r_l$ and $r_r$ of
$\mathfrak{g}_l := \text{Lie}_k ( G_l ( X ))$ and $\mathfrak{g}_r :=
\text{Lie}_k ( G_r ( X ))$ such that the Schouten algebraic brackets $[ r_l,
r_l ]$ and $[ r_r, r_r ]$ are equal as elements of $\mathfrak{g}_l \otimes_k k
\cong \mathfrak{g}_r \otimes_k k$ where $k$ is an algebraic closure of $k$.
Indeed, in this case, the bivector field $r_l^L - r_r^R$ is a Poisson bivector
on $X$, where $r_l^L$ (resp. $r_r^R$) is the left action of  $r_l$ (resp.
$r_r$) on $X$. Therefore, this problem is intimely linked with the problem of
classification of $r$-matrices in Lie algebras over a non-algebraically closed
field.

\begin{defin}
  A quantization of a $k$-Poisson torsor $( T_{\text{cl}}, m_{T_{\text{cl}}},
  1_{T_{\text{cl}}}, \mu_{T_{\text{cl}}}, \theta_{T_{\text{cl}}}, \{,
  \}_{\text{cl}} )$ is the data of a $k [[ \hbar ]]$-torsor $( T, m_T, 1_T,
  \mu_T$$, \theta_T )$ such that the $k [[ \hbar ]]$-algebra $( T, m_T, 1_T )$
  is a quantization of the $k$-Poisson algebra $( T_{\text{cl}},
  m_{T_{\text{cl}}}, 1_{T_{\text{cl}}}, \{, \}_{\text{cl}} )$ and such that
  $\mu_T = \mu_{T_{\text{cl}}}$ $( \text{mod} \hbar )$ and $\theta_T =
  \theta_{T_{\text{cl}}}$ (mod $\hbar )$.
\end{defin}

In the future, our goal is to classify {\em{all}} torsors corresponding to
certain classes of Hopf algebras. The problem seems to be difficult. Indeed,
when the ground field $k$ is not (necessarily) algebraically closed, affine
torsors of the algebraic geometry are classified by means of non-abelian
cohomology. On the other hand, when $k$ is algebraically closed, the study of
classical torsors in the category of Poisson manifolds {\cite{DS}} shows that
the Drinfeld twists should play a role in the classification.

\section{The results}

We present below our main results about torsors. As before, in all this
section, $k$ stands for a commutative field and $A$ for a commutative
$k$-algebra without any zero divisor.

\subsection{The Hopf algebras $H_l ( T )$ and $H_r ( T )$}

In the classical case, on any torsor, there are two groups acting simply and
transitively. In the non-commutative setting, we can also find two Hopf
algebras co-acting on a given torsor.

\begin{thm}
  \label{premtheo}Let $( T, m_T, 1_T, \mu_T, \theta_T )$ be an A-torsor. We
  denote by $H_l ( T, m_T, 1_T, \mu_T, \theta_T )$ or shortly $H_l ( T )$ or
  $H_l ( \mu_T )$ if there is no confusion, the subset of 
  $T \otimes_A T^{\text{op}}$ defined as
  $\left\{ x \in T \otimes_A
  T^{\text{op}} / ( \text{Id}_T \otimes \text{Id}_{T^{\text{op}}}\otimes
  \theta_T \otimes \text{Id}_{T^{\text{op}}} ) \circ ( \mu_T \otimes
  \text{Id}_{T^{\text{op}}} ) ( x ) = ( \text{Id}_T \otimes \mu_T^{\text{op}}
  ) ( x ) \right\}$. Then,
  \begin{enumerate}
    \item if $x \in H_l ( T ), m_T ( x )$ is a scalar denoted by
    $\varepsilon_{H_l ( T )} ( x ) 1_T$;
    
    \item if $x \in H_l ( T ), \Delta_{H_l ( T )} ( x ) := ( \mu_T
    \otimes \text{Id}_{T^{\text{op}}} ) ( x ) \in H_l ( T ) \otimes_A H_l ( T
    )$;
    
    \item if $x \in H_l ( T ),S_{H_l ( T )} ( x ) := \tau_{( 12 )}
    \circ ( \text{Id}_T \otimes \theta_T ) ( x ) \in H_l ( T )$;
    
    \item By defining $m_{H_l ( T )}$ as the restriction of $m_T \otimes
    m_T^{\text{op}}$ to $H_l ( T )$ and $\eta_{H_l ( T )} : A \rightarrow H_l
    ( T )$ as given by $\eta_{H_l ( T )} ( 1 ) = 1_T \otimes 1_T$, then the
    following sextuple 
    $( H_l ( T ), m_{H_l ( T )}, \Delta_{H_l ( T )}, \eta_{H_l ( T
    )}, \varepsilon_{H_l ( T )}, S_{H_l ( T )} )$ is an $A$-Hopf-algebra.
    
    \item $\text{Im} \mu_T \subset H_l ( T ) \otimes_A T$ and $\mu_T : T
    \longrightarrow H_l ( T ) \otimes_A T$ gives to $T$ a left $H_l ( T
    )$--comodule-algebra structure.
  \end{enumerate}
\end{thm}

  Likewise, we denote by $H_r ( T, m_T, 1_T, \mu_T, \theta_T )$ or shortly $H_r ( T )$ or
  $H_r ( \mu_T )$ if there is no confusion, the subset of
  $T^{\text{op}} \otimes_A T$ defined as 
  $$H_r ( T )
  := \left\{ x \in T^{\text{op}} \otimes_A T / ( \text{Id}_{T^{\text{op}}}
  \otimes \theta_T \otimes \text{Id}_{T^{\text{op}}} \otimes \text{Id}_{T^{}} )
  \circ ( \text{Id}_{T^{\text{op}}} \otimes \mu_T ) ( x ) = ( \mu_T^{\text{op}}
  \otimes \text{Id}_T ) ( x ) \right\}.$$ 
  Then, $H_r(T)$ can be equipped with a natural structure of $A$-Hopf algebra, 
  $\text{Im} \mu_T \subset T \otimes_A H_r ( T )$ and the map 
  $\mu_T : T \longrightarrow T \otimes_A H_r
  ( T )$ defines a right $H_r ( T )$-comodule-algebra structure on $T$.
Moreover, the two co-actions of $H_l(T)$ and $H_r(T)$ on $T$ commute.

\begin{example}
  \label{isotriv}Let $( H, m_H, \Delta_H, \eta_H, \varepsilon_H, S_H )$ be a
  Hopf algebra equipped with its trivial torsor structure $( H, m_H, 1_H,
  \mu_H, \theta_H )$. Then, $H_l ( H ) = ( \text{Id}_H \otimes S_H ) \circ
  \Delta ( H ), H_r ( H ) = ( S_H \otimes \text{Id}_H ) \circ \Delta_H ( H )$
  and $i_{l, H} := ( \text{Id}_H \otimes \varepsilon_H )$ (resp. $i_{r,
  H} := ( \varepsilon_H \otimes \text{Id}_H )$) establishes a Hopf
  algebra isomorphism between $H_l ( H )$ and $H$ (resp. $H_r ( H )$ and $H$).
\end{example}

\begin{example}
  \label{abelh}If $K / k$ is a Galois extension of a field $k$ equipped with
  its torsor structure seen in \ref{abel}, then $H_l ( T ) \cong ( k
  G^{\text{op}} )^{\ast}$ and $H_r ( T ) \cong ( k G )^{\ast}$ as Hopf
  algebras.
  
  \begin{example}
    \label{zeze}Let $A_{\alpha, \beta}$ be the non-commutative algebra-torsor
    considered in Example \ref{examp}. Then, $H_l ( A_{\alpha, \beta} ) = H_r
    ( A_{\alpha, \beta} )$ is generated by the elements $x \otimes x^{- 1}$
    and $y \otimes y^{- 1}$ and thus is isomorphic to the algebra of functions
    on $\mathbb{Z} / n \mathbb{Z} \times \mathbb{Z} / n \mathbb{Z} .$
  \end{example}
\end{example}

 \begin{note}
   \label{isoopp}If $( T, m_T, 1_T, \mu_T, \theta_T )$ is an $A$-torsor, then
   the map $( \theta_T \otimes \text{Id}_T )$ (resp. $( \text{Id}_T \otimes
   \theta_T )$) is a Hopf algebra isomorphism from $H_l ( \mu_T )$ to $H_r (
   \mu_T^{\text{op}} )$ (resp. $H_r ( \mu_T )$ to $H_l ( \mu_T^{\text{op}} )$).
 \end{note}
  
  \begin{note}
    \label{meme}If $( T, m_T, 1_T, \mu_T, \theta_T )$ is an $A$-torsor endowed
    with a commutative law, then $H_l ( T ) = H_r ( T )$ (in $T \otimes_A T$).
  \end{note}
  
  \begin{note}
    Let $( T, m_T, 1_T, \mu_T, \theta_T )$ be an $A$-torsor and $\varepsilon :
    T \longrightarrow A$ a character. If $\dim_k T < \infty$ or $\varepsilon
    \circ \theta_T = \varepsilon$, then $T$ is isomorphic {\em{as an
    algebra}} to its left or right Hopf algebra with an identification between
    $\varepsilon$ and the co-unity of the Hopf algebra.
  \end{note}
  
  \subsection{\label{toto}Torsors and Hopf-Galois extensions}
  
  The notion of Hopf-Galois extension was introduced by Chase and Sweedler in
  1969 in the commutative case {\cite{CS}} and later by Kreimer and Takeuchi
  in the non-commutative case {\cite{KT}} as a generalization of the Galois
  theory where the Galois groups are replaced by Hopf algebras. Let $k$ be a
  (commutative) field and $H$ an Hopf algebra. The axioms for a left
  $H$-Galois extension $A$ of $k$ are the ones we would get by taking formula
  (\ref{deftorgeoalg}) as a reference for defining torsors over Hopf algebras
  {\cite{Hjs}}. By definition, a left $H$-Galois extension $A$ of $k$ is a
  $H$-left comodule-algebra $A$ given by a morphism $\Delta_A : A
  \longrightarrow H \otimes A$ such that 1) $\left\{ x \in A / \Delta_A ( x )
  = 1 \otimes x \right\} = k$ and 2) the natural map $\text{can} : A \otimes A
  \longrightarrow H \otimes A$ given by $\text{can} = ( \text{Id}_H \otimes
  m_A ) \circ ( \Delta_A \otimes \text{Id}_A )$ is a {\em{linear}}
  isomorphism. In the classical case, if $X$ is a torsor, then the natural map
  \[ \begin{array}{ccccc}
       G_l ( X ) \times X &  & \longrightarrow &  & X \times X\\
       ( g, x ) &  & \longmapsto &  & ( g . x, x )
     \end{array} \]
  is bijective. The following proposition generalizes this fact to the
  non-commutative case.

  \begin{prop}
    \label{hopfgalprop}Let $( T, m_T, 1_T, \mu_T, \theta_T )$ be a $k$-torsor.
    Then, $T$ is a left $H_l ( T )$--Galois extension of $k$ and a right $H_r
    ( T )$-Galois-extension of $k$. In particular, in the finite dimensional
    case, we have $\dim_k T = \dim_k H_l ( T ) = \dim_k H_r ( T )$.
  \end{prop}

\subsection{Composition of torsors}

The $A$-torsors form a category. We can naturally define the notion of torsors
morphism, sub-torsors, quotient torsors and tensor product of two torsors. If
$( T_i, m_{T_i}, 1_{T_i}, \mu_{T_i}, \theta_{T_i} ), i = 1, 2$ are two torsors
with a Hopf algebra isomorphism between $H_l ( T_1 )$ and $H_r ( T_2 )$, then
the following theorem shows that we can compose the two torsors to get a third
torsor whose Hopf algebra co-acting from the left (resp. from the right) is
isomorphic to $H_l ( T_1 )$ (resp. $H_r ( T_2 )$).

\begin{thm}
  \label{zorro}\label{compo}Let $( T_i, m_{T_i}, 1_{T_i},
  \mu_{T_i}, \theta_{T_i} ), i = 1, 2$ be two $A$-torsors where $A$ is as
  above. Let us assume that there is an isomorphism $\Phi : H_r ( T_1 )
  \longrightarrow H_l ( T_2 )$ of $A$-Hopf algebra and set 
  $$T_{\Phi} :=
  T_1 \otimes_{\Phi} T_2 := \left\{ x \in T_1 \otimes_A T_2 / (
  \text{Id}_{T_1} \otimes \Phi \otimes \text{Id}_{T_2} ) \circ ( \mu_{T_1}
  \otimes \text{Id}_{T_2} ) ( x ) = ( \text{Id}_{T_1} \otimes \mu_{T_2} ) ( x
  ) \right\}.$$ Let $m_{\Phi}$ and $\theta_{\Phi}$ be the restrictions of
  $m_{T_1} \otimes m_{T_2}$ and $\theta_{T_1} \otimes \theta_{T_2}$ to
  $T_{\Phi}$ and let $\mu_{\Phi}$ be the map defined on $T_{\Phi} \subset T_1
  \otimes_A T_2$ with the help of the generalized Sweedler notations by :
  \begin{equation}
    \mu_{\Phi} ( x_i \otimes y_i ) = \tau_{( 34 )} ( x_i^{( 1 )} \otimes \Phi
    ( x_i^{( 2 )} \otimes x_i^{( 3 )} ) \otimes x_i^{( 4 )} \otimes x_i^{( 5
    )} \otimes y_i )
  \end{equation}
  where $\tau_{( 34 )} : T_1 \otimes_A T_2 \otimes_A T_2^{\text{op}} \otimes_A
  T_1 \otimes_A T_1 \otimes_A T_2 \longrightarrow T_1 \otimes_A T_2 \otimes_A
  T_1 \otimes_A T_2^{\text{op}} \otimes_A T_1 \otimes_A T_2$ denotes the
  permutation morphism of the third and fourth factors, $H_l ( T_2 )$ being
  imbedded in $T_2 \otimes_A T_2^{\text{op}}$. Then,
  \begin{enumerate}
    \item $( T_{\Phi}, m_{\Phi}, 1_{T_1} \otimes 1_{T_2} )$ is an $A$-algebra,
    $\text{Im} \mu_{\Phi} \subset T_{\Phi} \otimes T_{\Phi}^{\text{op}}
    \otimes T_{\Phi}, \text{Im} \theta_{\Phi} \subset T_{\Phi}$ and $(
    T_{\Phi}, m_{\Phi}, 1_{T_1} \otimes 1_{T_2}$$, \mu_{\Phi}, \theta_{\Phi}
    )$ is an $A$-torsor.
    
    \item The restriction of the map $\tau_{( 34 )} \circ ( \text{Id}_{T_1}
    \otimes \Phi \otimes \text{Id}_{T_1^{\text{op}}} ) \circ ( \mu_{T_1}
    \otimes \text{Id}_{T_1^{\text{op}}} )$ to $H_l ( T_1 ) \subset T_1
    \otimes_A T_1^{\text{op}}$ gives rise to an $A$-Hopf algebra isomorphism
    between $H_l ( T_1 )$ and $H_l ( T_{\Phi} )$, whose inverse map is $(
    \text{Id}_{T_1} \otimes \varepsilon_{H_l ( T_2 )} \otimes
    \text{Id}_{T_1^{\text{op}}}$$) \circ \tau_{( 34 )}$.
    
    \item The restriction of the map $\tau_{( 12 )} \circ (
    \text{Id}_{T_2^{\text{op}}} \otimes \Phi^{- 1} \otimes \text{Id}_{T_2} )
    \circ ( \text{Id}_{T_2^{\text{op}}} \otimes \mu_{T_2} )$ to $H_r ( T_2 )
    \subset T_2^{\text{op}} \otimes_A T_2$ gives rise to an $A$-Hopf algebra
    isomorphism between $H_r ( T_2 )$ and $H_r ( T_{\Phi} )$, whose inverse
    map is $( \text{Id}_{T_2^{\text{op}}} \otimes \varepsilon_{H_r ( T_1 )}
    \otimes \text{Id}_{T_2} ) \circ \tau_{( 12 )}$.
  \end{enumerate}
\end{thm}

\subsection{The group $\text{Tor} ( H )$}

If $f : ( T_1, m_{T_1}, 1_{T_1}, \mu_{T_1}, \theta_{T_1} ) \longrightarrow (
T_2, m_{T_2}, 1_{T_2}, \mu_{T_2}, \theta_{T_2} )$ is an $A$-torsor morphism,
then it can be shown that $( f \otimes f^{\text{op}} ) ( H_l ( T_1 )) \subset
H_l ( T_2 )$ and that $f_l := ( f \otimes f^{\text{op}} )_{| H_l ( T_1 )}
: H_{l ( T_1 )} \longrightarrow H_l ( T_2 )$ is an $A$-Hopf algebra morphism.
Likewise, we define a Hopf morphism $f_r : H_r ( T_1 ) \longrightarrow H_r (
T_2 )$. Consequently, if $H$ and $H'$ are two $A$-Hopf algebras, on the set
$\widehat{\text{Tor}} ( H, H' )$ of the septuples $( T_{}, m_{T_{}}, 1_{T_{}},
\mu_{T_{}}, \theta_{T_{}}, i_{l, T_{}}, i_{r, T_{}} )$where $( T_{}, m_{T_{}},
1_{T_{}}, \mu_{T_{}}, \theta_{T_{}} )$ is an $A$-torsor and where $i_{l, T_{}}
: H_l ( A ) \longrightarrow H$ and $i_{r, T} : H_r ( A ) \longrightarrow H'$
are two $A$-Hopf algebra isomorphisms, we can define a relation $\sim_{H, H'}$
by $( T_1, m_{T_1}, 1_{T_1}, \mu_{T_1}, \theta_{T_1}, i_{l, T_1}, i_{r, T_1} )
\sim_{H, H'} ( T_2, m_{T_2}, 1_{T_2}, \mu_{T_2}, \theta_{T_2}, i_{l, T_2},
i_{r, T_2} )$ if and only if there exists an $A$-torsor {\em{isomorphism}}
$f : ( T_1, m_{T_1}, 1_{T_1}, \mu_{T_1}, \theta_{T_1} ) \longrightarrow ( T_2,
m_{T_2}, 1_{T_2}, \mu_{T_2}, \theta_{T_2} )$ such that $i_{l, T_1} = i_{l,
T_2} \circ f_l $ and $i_{r, T_1} = i_{r, T_2} \circ f_r$. This relation is an
{\em{equivalence}} relation. The quotient set is denoted by $\text{Tor} ( H,
H' )$. Moreover, if $H$, $H'$, $H''$ are three $A$-Hopf algebras, we have a
natural map :
\begin{equation}
  \begin{array}{ccc}
    \label{defcomp} \widehat{\text{Tor}_{}} ( H, H' ) \times
    \widehat{\text{Tor}_{}} ( H', H'' ) & \longrightarrow &
    \widehat{\text{Tor}} ( H, H'' )\\
    ( T_1,\ldots, i_{r, T_1}
    )\times ( T_2,\ldots
    i_{r, T_2} ) & \longmapsto & ( T,\ldots
    i_{r, T} )
  \end{array} 
\end{equation}
where $T := T_1 \otimes_{\Phi} T_2$ is by Theorem \ref{zorro} an
$A$-torsor and $\Phi := i_{l, T_2}^{- 1} \circ i_{l, T_1}$, $i_{l, T}
:= i_{l, T_1} \circ( \text{Id}_{T_1} \otimes \varepsilon_{H_l ( T_2 )}
\otimes \text{Id}_{T_1^{\text{op}}}) \circ \tau_{( 34 )}$ and $i_{r, T}
:= i_{r, T_2} \circ ( \text{Id}_{T_2^{\text{op}}} \otimes
\varepsilon_{H_r ( T_1 )} \otimes \text{Id}_{T_2} ) \circ \tau_{( 12 )}$. It
can be shown that the map defined in (\ref{defcomp}) is associative and
compatible with the equivalence relation $\sim_{H, H'}$. Thus, we define a
composition law on $\text{Tor} ( H ) := \text{Tor} ( H, H )$ which is
denoted by $\ast$.

\begin{thm}
  \label{invtor}Let $( H, m_H, \Delta_H, \eta_H, \varepsilon_H, S_H )_{}$ be a
  Hopf algebra over a commutative field $k$. The set $\text{Tor} ( H )$
  equipped with the law $\ast$ is a group whose unit element is the class of the
  trivial torsor $( H, m_H, 1_H, \mu_H, \theta_H, i_{l, H}, i_{r, H} )$ (the
  notations are the same as those of Examples \ref{tortriv} and
  \ref{isotriv}). The inverse of the class of $( T, m_T, 1_T, \mu_T,
  \theta_T, i_{l, T}, i_{r, T} )$ is equal to the class of the opposite torsor
  of $T$ : 
  $( T^{\text{op}}, m_T^{\text{op}}, 1_T, \mu^{\text{op}}_T, \theta_T,
  i_{l, T^{\text{op}}}, i_{r, T^{\text{op}}} )$ with $i_{l, T^{\text{op}}}
  := i_{r, T} \circ ( \text{Id}_T \otimes \theta_T )^{- 1}$ and $i_{r,
  T^{\text{op}}} := i_{l, T} \circ ( \theta_T \otimes \text{Id}_T )^{-
  1}$ (see Remark \ref{opp} and Note \ref{isoopp}).
\end{thm}

The group $\text{Tor} ( H )$ is called the torsor invariant of the $A$-Hopf
algebra $H$. As we saw with the study of Example \ref{examp} together with
Example \ref{zeze}, this group is far from trivial. Given the link with the
Hopf-Galois extensions theory, it should be seen as a subgroup of
$\text{Bigal} ( H )$ introduced earlier by Schauenburg {\cite{Sch}} and also
as a generalization of the Harisson group in Galois theory.

\subsection{Cotorsors ``à la Parmentier''}

We show that Parmentier's formalism which allows to quantize affine Poisson
groups {\cite{DS}} can be subsumed in our theory of torsors or cotorsors whose
definition is given below.

\begin{defin}
  An $A$-cotorsor is a quintuple $( C, \Delta_C, \varepsilon_C, \nu_C,
  \theta_C )$ where $( C, \Delta_C, \varepsilon_C )$ is an $A$-coalgebra,
  $\nu_C : C \otimes_A C^{\text{cop}} \otimes_A C \longrightarrow C$ an
  $A$-coalgebra morphism and $\theta_C : C \longrightarrow C$ a coalgebra
  automorphism satisfying the following axioms :
  \begin{eqnarray}
    \nu_C \circ ( \Delta_C \otimes \text{Id}_C ) & = & \varepsilon_C \otimes
    \text{Id}_C\\ 
    \nu_C \circ ( \text{Id}_C \otimes \Delta_C ) & = & \text{Id}_C \otimes
    \varepsilon_C\\ 
    \nu_C \circ ( \nu_C \otimes \text{Id}_{C^{\text{cop}}} \otimes \text{Id}_C
    ) & = & \nu_C \circ ( \text{Id}_C \otimes \text{Id}_{C^{\text{cop}}}
    \otimes \nu_C )\\ 
    \nu_C \circ ( \nu_C \otimes \text{Id}_{C^{\text{cop}}} \otimes \text{Id}_C
    ) \circ  \theta^{(3)}_C& = & \nu_C \circ
    ( \text{Id}_C \otimes \nu_C^{\text{op}} \otimes \text{Id}_C )\\ 
    \nu_C \circ ( \theta_C \otimes \theta_C \otimes \theta_C ) & = & \nu_C
    \circ \theta_C 
  \end{eqnarray}
  with $\nu_C^{\text{op}} := \nu_C \circ \tau_{( 13 )}$ and
$\theta^{(3)}_C:= ( \text{Id}_C \otimes \text{Id}_{C^{\text{cop}}} \otimes \theta_C
    \otimes \text{Id}_{C^{\text{cop}}} \otimes \text{Id}_C )$.
\end{defin}

The cotorsors theory can be developed in the same way as the torsors one. In
particular, any cotorsor can be equipped with two coalgebra-module structures
over Hopf algebras and the two actions commute. If $( T, m_T, 1_T, \mu_T,
\theta_T )$ is an $A$-torsor, then $( T^{\ast}, m_T^{\ast}, \eta_T^{\ast},
\mu_T^{\ast}, \theta_T^{\ast} )$ is an $A$-cotorsor, where $\eta_T : A
\longrightarrow T$ is defined by $\eta_T ( 1 ) = 1_T$. The converse is true in
the finite dimentional case.

Now, let $( H, m_H, \Delta_H, \eta_H, \varepsilon_H, S_H )$ be an $A$-Hopf
algebra and let $F \in H \otimes_A H$ be a Drinfeld twist {\cite{Dri}}, i.e.,
an element $F$ satisfying the equations :
\begin{equation}
  ( F \otimes 1 ) ( \Delta \otimes \text{Id}_H ) ( F ) = ( 1 \otimes F ) (
  \text{Id}_H \otimes \Delta_H ) ( F ) \label{twist}
\end{equation}
\begin{equation}
  ( \varepsilon_H \otimes \text{Id}_H ) ( F ) = ( \text{Id}_H \otimes
  \varepsilon_H ) ( F ) = 1 \label{twistbis}
\end{equation}

Then, it is known that $u_F := m_H \circ ( \text{Id}_H \otimes S_H ) ( F
)$ is an invertible element of $H$ whose inverse is $u_F^{- 1} = m_H \circ (
S_H \otimes \text{Id}_H ) ( F^{- 1} )$ and that if we set $\Delta_F := F
\Delta_H F^{- 1}$ and $S_F := u_F S_H u_F^{- 1}$, then the sextuple $( H,
m_H, \Delta_F, \eta_H, \varepsilon_H, S_F )$ denoted shortly by $H_F$ is an
$A$-Hopf algebra. Moreover, if $H = U_{\hbar} ( \mathfrak{g} )$ is a QUE algebra, 
$\delta : = \lim_{\hbar \rightarrow 0} \hbar^{-
1} ( \Delta_H - \Delta_H^{\text{op}} )$ and $f := \lim_{\hbar \rightarrow
0} \hbar^{- 1} ( F - 1 \otimes 1 )$, then the triple $( H, \Delta_H F^{- 1},
\varepsilon_H )$ is a $k$-coalgebra which is a quantization of the affine
Poisson group given by the triple $( \mathfrak{g}, \delta, f )$ {\cite{Par}}.

\begin{thm}
  \label{parmi}Let $( H, m_H, \Delta_H, \eta_H, \varepsilon_H, S_H )$ be an
  $A$-Hopf algebra and let $F \in H \otimes_A H$ be a Drinfeld twist. Set $(
  C, \Delta_C, \varepsilon_C ) := ( H, \Delta_H F^{- 1}, \varepsilon_H )$
  and let $\theta_C$ be the $A$-linear map defined on $C$ by $\forall x \in
  C,\theta_C ( x ) = S_H^2 ( x ) S_H ( u_F ) u_F^{- 1}$ with $u_F :=
  m_H \circ ( \text{Id}_H \circ S_H ) ( F )$ and
  \begin{equation}
    \begin{array}{cccc}
      \nu_C : & C \otimes_A C^{\text{cop}} \otimes_A C & \longrightarrow & C
      \label{nu}\\
      & x \otimes y \otimes z & \longmapsto & x u_F S_H ( y ) z
    \end{array}
  \end{equation}
  Then $\theta_C$ is an $A$-coalgebra automorphism for the $A$-coalgebra $( C,
  \Delta_C, \varepsilon_C ), \nu_C$ is a $A$-coalgebra morphism, and the
  quintuple $( C, \Delta_C, \varepsilon_C, \nu_C, \theta_C )$ is an
  $A$-cotorsor whose Hopf algebra acting from the left (resp. from the right)
  on $C$ is isomorphic to $H$ (resp. $H_F$).
\end{thm}

In particular, if we denote by $i_{l, C} $ and $i_{r, C}$ two isomorphisms
from $H$ and $H_F$ to $H_l ( C )$ and $H_r ( C )$, we see that, in the case
where $A = k$ and $C$ is a finite dimensional $k$-vector space, then the
septuple $( C^{\ast}, \Delta_C^{\ast}, \varepsilon_C^{\ast} ( 1 ),
\nu_C^{\ast}, \theta_C^{\ast}, \text{$i_{l, C}^{\ast}$$^{},$$i_{r, C}^{\ast}$}
) \in \widehat{\text{Tor}_{}} ( H, H_F )$. Moreover, if $F$ commutes with the
comultiplication $\Delta_H$, then $H_F = H$. Therefore, we see that
$\widehat{\text{Tor}} ( H )$ contains a subgroup isomorphic to $\text{Aut} ( H
)^2 \times \{ F \in ( H \otimes H )^{\times} / F$ is a Drinfeld twist and $F$
commutes with $\text{Im} ( \Delta_H ) \}$. On the other hand, according to
Note \ref{meme}, to get a torsor which is not of the Parmentier type, it is
enough to consider torsors endowed with a commutative law whose underlying
algebra does not possess any character. This is the case of torsors of Example
\ref{quadra} or Example \ref{examp}. This last example gives also torsors
which are neither ``Parmentier type'' torsors nor torsors arising from
algebraic geometry.

\section{Open problems}

\begin{enumerate}
  \item Find a {\em{non-trivial}} Poisson torsor structure on a commutative
  algebra without any character and try to quantize it. Classify all Poisson torsors.
  
  \item Find $\text{Tor} ( H )$ in simple cases.
  
  \item Study $\text{Tor} ( H )$ when the Hopf algebra $H$ is a quantized
  enveloping algebra.
  
  \item In particular, if $( \mathfrak{g}, \delta )$ a $k$-Lie bialgebra, if
  $f$ is a classical Drinfeld twist and if $( U_{\hbar} ( \mathfrak{g} ),
  \Delta )$ is a quantization of $( \mathfrak{g}, \delta )$, is it true or not
  that every quantization of the Poisson torsor given by the triple $(
  \mathfrak{g}, \delta, f )$ is a ``Parmentier'' type torsor, that is to say
  given by a quantization $F$ of $f$ ?
  
  \item Study the case where $H = \mathbb{C} \Gamma$is equipped with its
  natural Hopf algebra structure, with $\Gamma$ a finite group, by using some
  results due to Movshev {\cite{M}} or to Etingof-Gelaki {\cite{EG}} about
  twists in this Hopf algebra.
  
  \item Find and classify all ``low dimensional'' torsors.
  
  \item Generally speaking, classify all torsors and cotorsors for a large 
        class of Hopf algebras.
  
  \item A classical torsor together with its opposite torsors can be seen as a 
        groupoïd on a basis with two elements. Are our axioms compatible with 
        the notion of a quantum groupoïd {\cite{Mal}} ?
  
  \item Is it possible to find an Hopf-Galois extension which is not a torsor 
        in our sense ?
  \item Establih a link with the Galois cohomology which classifies torsors in
  algebraic geometry.
\end{enumerate}

\small{
{}
}

\end{document}